\documentclass[fleqn,11pt]{article}
\usepackage[numbers,sort&compress]{natbib}
\usepackage[footnotesep=0.4in]{geometry}
\usepackage{graphicx}
\usepackage{amsmath}
\usepackage{bbding}
\usepackage{setspace}
\usepackage{pifont}
\usepackage{multirow}
\usepackage{amsfonts}
\usepackage{caption2}
\usepackage{float}
\usepackage{graphics}
\usepackage{subfigure}
\usepackage{float}
\usepackage{array}
\usepackage{booktabs}
\makeatletter
\newcommand\figcaption{\def\@captype{\textbf{figure}}\caption}
\newcommand\tabcaption{\def\@captype{table}\caption}
\begin{document}
\begin{spacing}{1.0}
\date{}
\title{Numerical simulation for the motions of nonautonomous solitary waves of a variable-coefficient forced Burgers equation via the lattice Boltzmann method}
\author{Qing-Feng Guan$^{1,\thanks{Corresponding
author: vguan07@163.com}}$, Wei-Qin Chen$^2$,  Ying Li$^3$
\\{\em \footnotesize School of Mathematics, Shandong University, Jinan, China, 250100}
\\{\em \footnotesize Department of Electrical, Computer, and Systems Engineering, Rensselaer Polytechnic Institute, NY, USA, 12180}
\\{\em \footnotesize School of Computer and Communication Engineering, Changsha University of Science and Technology , Changsha, China, 410114}
}
\maketitle

\vspace{-10mm}
\begin{abstract}
\vspace{-2mm}
\renewcommand{\raggedright}{\leftskip=0pt \rightskip=0pt plus 0cm}
The lattice Boltzmann method (LBM) for the variable-coefficient forced Burgers equation (vc-FBE) is studied by choosing the equilibrium distribution and compensatory functions properly. In our model, the vc-FBE is correctly recovered via the Chapman-Enskog analysis. We numerically investigate the dynamic characteristics of solitons caused by the dispersive and external-force terms. Four numerical examples are given, which align well with the theoretical solutions. Our research proves that LBM is a satisfactory and efficient method for nonlinear evolution equations with variable coefficients.


\end{abstract}


\vspace{-6mm}
\section{Introduction}\label{S1}
\vspace{-1mm}

The lattice Boltzmann method (LBM) has drawn broad attention in simulations of fluid flows in different mathematical and physical fields~\cite{1}. As one of the most popular models via the LBM, the lattice Bhatnagar-Gross-Krook one has been applied in various fields successfully~\cite{2}. Compared with other traditional simulation approaches (for instance, the finite difference~\cite{3}, Galerkin~\cite{4}, Chebyshev spectral~\cite{5}, finite element~\cite{6}, heat balance integral~\cite{7}, variational iteration~\cite{8} and decomposition methods~\cite{9}), the LBM has such advantages as the parallel calculation, algorithmic simplicity and easy-to-handle of marginal conditions. Based on the successive Boltzmann equation, the macroscopic variable is derived easily and the macroscopic equation is also recovered correctly under the assistance of Chapman-Enskog expansion~\cite{10}.

Recently, the LBM has been employed in certain nonlinear evolution equations (NLEEs), including the Burgers~\cite{11}, Korteweg-de Vries and Burgers~\cite{12}, improved Korteweg-de Vries~\cite{13} and convection-diffusion equations~\cite{14}, to name a few. Nevertheless, such studies mainly focused on the constant-coefficient cases but little involved the variable-coefficient ones, especially for the nonlinear models with variable external-force term. External-force coefficient can 
depict the time evolution of varying interface contour in different reaction-diffusion models~\cite{15}. Indeed, Ref.~\cite{16} shows that the nonautonomous solitons can show more dynamic characteristics, which do not appear in the constant-coefficient NLEEs~\cite{17,18,19,20}.

In this work, we work over the variable-coefficient forced Burgers equation (vc-FBE) as follows~\cite{21}
\setlength\abovedisplayskip{0.5cm}
 \setlength\belowdisplayskip{0.5cm}
\vspace{-0.1mm} 
\begin{equation}\label{Eq1}
\begin{aligned}
&u_{t}+a(x,t) u u_{x}+ b(x,t) u_{xx}-m(x,t)=0,
\\
\end{aligned}
\end{equation}
\vspace{-0.1mm}
where the function $u(x,t)$ represents the fluid velocity field, and $a(x,t)$, $b(x,t)$ and $m(x,t)$ denote the nonlinear, dispersive and external-force coefficients, respectively. Xu $et$ $al$. have presented the multi-soliton solutions of Eq.~(\ref{Eq1}) through the generalized Hopf-Cole transform~\cite{15}. Gao $et$ $al$.~\cite{21} have achieved several kinds of exact soliton-like solutions of Eq.~(\ref{Eq1}). When $a(x,t)$ = 1, $b(x,t) = -\frac{1}{Re}$ and $m(x,t)$ = 0, Eq.~(\ref{Eq1}) can be 
transformed into  the Burgers equation~\cite{22}
\vspace{-2mm}
\begin{equation}\label{Eq2}
\setlength\abovedisplayskip{0.5cm}
\setlength\belowdisplayskip{0.5cm}
\begin{aligned}
&u_{t}+ u u_{x}- \frac{u_{xx}}{Re}=0,
\\
\end{aligned}
\end{equation}
\vspace{-0.2mm}
where $Re$ stands for the Reynolds number describing the coefficient of kinematic viscosity $v$ ($v= \frac{1}{Re}$). Refs.~\cite{23,24,25} have numerically investigated the Burgers equation through the LBM. However, to the best of our knowledge, few studies focused on the Burgers equation with the variable and external-force coefficients based on the LBM.

In this work, by using the LBM, we will show the dynamics of the nonautonomous solitons described by Eq.~(\ref{Eq1}). The rest part of this article is listed as below: in Section~\ref{S2}, the lattice Boltzmann model will be proposed; in Section~\ref{S3}, the accurateness and durability of the proposed model will be verified by four numerical experiments; our consequences will be given in Section~\ref{S4}.

\vspace{-2mm}
\section{The lattice Boltzmann model}\label{S2}
\vspace{-1mm}

In this section, we choose the D1Q3 velocity model in which the discrete velocities can be written as: $\{c_0, c_1, c_2\}=\{0, c, -c\}$.
$c$ is a constant and $c=\Delta x/\Delta t$, in which $\Delta t$ and $\Delta x$ are the lattice time and space steps, singly. The evolution equation described by the density distribution function is shown below.
\vspace{-2mm}
\begin{equation}\label{Eq3}
\setlength\abovedisplayskip{0.3cm}
\setlength\belowdisplayskip{0.3cm}
\begin{aligned}
&f_{i}(x+c_{i}\Delta t ,t+\Delta t )=f_{i}(x,t)-\frac{1}{\tau}(f_{i}(x,t)-f_{i}^{eq}(x,t))+ h_{i}(x,t) \Delta t +\frac{ m(x,t)}{3} \Delta t,
\end{aligned}
\end{equation}
\vspace{-1mm}
in which $f_{i}(x,t)$ and $f_{i}^{eq}(x,t)$ are defined as the density distribution and equilibrium distribution functions and $\tau$ is the dimensionless relaxation time. $h_{i}(x,t)$ is a compensatory function added in Eq.~(\ref{Eq3})
\vspace{-2mm}
\begin{equation}\label{Eq4}
\begin{aligned}
&h_{i}(x,t)=\lambda_{i}u^2.
\end{aligned}
\end{equation}
\vspace{-5mm}

To generate the equilibrium distribution function, Taylor expansion is applied to Eq.~(\ref{Eq3}) with the terms expanded up to \rm{O}$(\Delta t^3)$

\vspace{-2mm}
\begin{equation}\label{Eq5}
\setlength\abovedisplayskip{0.3cm}
\setlength\belowdisplayskip{0.3cm}
\begin{aligned}
&\Delta t\frac{\partial f_{i}(x,t)}{\partial t}+\Delta t\frac{\partial(c_{i}f_{i}(x,t))}{\partial x}+\frac{1}{2}\Delta t^2\frac{\partial^2 f_{i}(x,t)}{\partial t^2}+\Delta t\frac{\partial^2(c_{i}f_{i}(x,t))}{\partial x\partial t}+\frac{1}{2}\Delta t^2\frac{\partial^2 (c_{i}^2f_{i}(x,t))}{\partial x^2}
\vspace{1ex}\\&
=-\frac{1}{\tau}(f_{i}(x,t)-f_{i}^{eq}(x,t))+\Delta t h_{i}(x,t)+\frac{\Delta t m(x,t)}{3}+{\rm O}(\Delta t^3).&
\end{aligned}
\end{equation}
\vspace{-1mm}
Then, we give the following expressions via the Chapman-Enskog expansion
\vspace{-1mm}
\begin{equation}\label{Eq6}
\setlength\abovedisplayskip{0.3cm}
\setlength\belowdisplayskip{0.3cm}
\begin{aligned}
&f_{i}(x,t)=\sum_{n=0}^{\infty}\epsilon^nf_{i}^{(n)}(x,t)=\epsilon^0 f_i^{\left( 0 \right)}(x,t)+\epsilon^1 f_i^{\left( 1 \right)}(x,t)+\epsilon^2 f_i^{\left( 2 \right)}(x,t)+... ~ ,
\\
&\frac{\partial}{\partial t}=\epsilon\frac{\partial}{\partial t_1},
m(x,t)=\epsilon m_{1}(x,t).
\end{aligned}
\end{equation}
Let $\epsilon$=$\Delta t$, by inserting Eq.~(\ref{Eq6}) into Eq.~(\ref{Eq5}), we can obtain a serial of differential equations
\\
\begin{equation}\label{Eq7}
\setlength\abovedisplayskip{0.3cm}
\setlength\belowdisplayskip{0.3cm}
\begin{aligned}
\rm{O} &(\Delta t^0): f_i^{\left( 0 \right)}(x,t) - f_i^{\left( {eq} \right)}(x,t)=0,
\end{aligned}
\end{equation}

\noindent
\begin{equation}\label{Eq8}
\setlength\abovedisplayskip{0.3cm}
\setlength\belowdisplayskip{0.3cm}
\begin{aligned}
\rm{O} &(\Delta t^1): \frac{\partial }{{\partial x}}\left( {{c_i}f_i^{\left( 0 \right)}(x,t)} \right) + \frac{1}{\tau }f_i^{\left( 1 \right)}(x,t) - h_{i}(x,t)=0,
\end{aligned}
\end{equation}

\noindent
\begin{equation}\label{Eq9}
\begin{aligned}
\rm{O} &(\Delta t^2): \frac{{\partial f_i^{\left( 0 \right)}(x,t)}}{{\partial {t_1}}} + \frac{\partial }{{\partial x}}\left( {{c_i}f_i^{\left( 1 \right)}(x,t)} \right) + \frac{1}{2}\frac{{{\partial ^2}}}{{\partial {x^2}}}\left( {c_i^2f_i^{\left( 0 \right)}(x,t)} \right) + \frac{1}{\tau }f_i^{\left( 2 \right)}(x,t) - \frac{{{m_1}\left( {x,t} \right)}}{3}=0.
\end{aligned}
\end{equation}

Similar to the general  lattice Boltzmann model, the macroscopic variable $u$ is given by
\begin{equation}\label{Eq10}
\begin{aligned}
&u=\sum\limits_i {f_i}(x,t).
\end{aligned}
\end{equation}
The equilibrium distribution function should suit the following constraint via the conservative law
\begin{equation}\label{Eq11}
\begin{aligned}
&u=\sum\limits_i {f_i^{\left( 0 \right)}(x,t)} = \sum\limits_i {f_i^{\left( {eq} \right)}(x,t)}.
\end{aligned}
\end{equation}
Based on Eq.~(\ref{Eq11}), we have
\begin{equation}\label{Eq12}
\begin{aligned}
&\sum\limits_i {f_i^{\left( n \right)}(x,t)}  = 0,\left( {n > 0} \right).
\end{aligned}
\end{equation}

Considering the structure of Eq.~(\ref{Eq1}), we have to take other constraints on the equilibrium distribution and compensatory functions
\begin{equation}\label{Eq13}
\begin{aligned}
&\sum\limits_i {c_if_i^{\left( 0 \right)}(x,t)}  = {F_1},
~\sum\limits_i {c_i^2f_i^{\left( 0 \right)}(x,t)}  = {F_2},
~\sum\limits_i {{h_i}(x,t) = {H_0}},
~\sum\limits_i {{c_i}{h_i}(x,t) = {H_1}}.
\end{aligned}
\end{equation}
Through Eq.~(\ref{Eq8}) and noting Eq.~(\ref{Eq13}), we obtain
\begin{equation}\label{Eq14}
\begin{aligned}
&\sum\limits_i {{c_i}f_i^{(1)}(x,t)} = \tau \sum\limits_i {{c_i}{h_i}(x,t) - \tau \frac{\partial }{{\partial x}}}  \sum\limits_i {c_i^2f_i^{(0)}(x,t)}  = \tau {H_1} - \tau \frac{\partial F_2}{\partial x}.
\end{aligned}
\end{equation}
By using the summation of Eq.~(\ref{Eq8}) over $i$, and considering Eqs.~(\ref{Eq11}) and (\ref{Eq13}), we get
\begin{equation}\label{Eq15}
\begin{aligned}
\frac{{\partial {F_1}}}{{\partial x}}={H_0}.
\end{aligned}
\end{equation}
Similarly, using the summation of Eq.~(\ref{Eq9}) over $i$, and noting Eqs.~(\ref{Eq11}), (\ref{Eq12}), (\ref{Eq13}) and (\ref{Eq14}), we have
\begin{equation}\label{Eq16}
\begin{aligned}
&\frac{{\partial u}}{{\partial {t_1}}} + (\frac{1}{2} - \tau )\frac{{{\partial ^2}{F_2}}}{{\partial {x^2}}}+ \tau \frac{{\partial {H_1}}}{{\partial x}} = {m_1}(x,t).
\end{aligned}
\end{equation}

In terms of the structure of $\Delta t$$\times$(\ref{Eq15})+$\Delta t^2$$\times$(\ref{Eq16}), Eq.~(\ref{Eq1}) can be converted into
\begin{equation}\label{Eq17}
\begin{aligned}
&\frac{{\partial u}}{{\partial t}} + \Delta t \tau \frac{{\partial {H_1}}}{{\partial x}} + \frac{{\partial {F_1}}}{{\partial x}} + (\frac{1}{2} - \tau ) \Delta t \frac{{{\partial ^2}{F_2}}}{{\partial {x^2}}} = m(x,t)+{H_0}.
\end{aligned}
\end{equation}
Comparing Eq.~(\ref{Eq17}) with Eq.~(\ref{Eq1}), we take
\begin{equation}\label{Eq18}
\begin{aligned}
&{F_1} = 0,~{F_2} = c^2 \eta u,
~{H_0} = 0,~{H_1} = c{\lambda}{u^2}.
\end{aligned}
\end{equation}
Applying Eq.~(\ref{Eq18}) to Eq.~(\ref{Eq17}), we finally have the following equation
\vspace{-1mm}
\begin{equation}\label{Eq19}
\begin{aligned}
&{u_t} + 2\tau \Delta t c \lambda u{u_x} + ( \frac{1}{2} - \tau)\Delta t {c^2} \eta u_{xx} = m(x,t).
\end{aligned}
\end{equation}
Observing the forms of Eqs.~(\ref{Eq19}) and (\ref{Eq1}), we get
\begin{equation}\label{Eq20}
\begin{aligned}
&\begin{array}{l}
a(x,t) = 2\tau \Delta t c \lambda,
~b(x,t) = ( \frac{1}{2} - \tau) \Delta t {c^2} \eta.
\end{array}
\end{aligned}
\end{equation}
Combining Eq.~(\ref{Eq10}) and the first two terms of Eq.~(\ref{Eq13}), the equilibrium distribution function is written as
$$
f_{i}^{(eq)}(x,t)=
\begin{cases}
~(1-\eta) u,&\text{$i$=0}, \\
\quad ~~ \frac{\eta u}{2}, &\text{$i$=1}, \\
\quad ~~ \frac{\eta u}{2}, &\text{$i$=2}.
\end{cases}
$$
Combining Eqs.~(\ref{Eq4}), (\ref{Eq18}) and the last two terms of Eq.~(\ref{Eq13}), the compensatory function $h_{i} (x,t)$ can be solved. To simplify the simulation without losing generality, we take only one case for $h_{i}(x,t)$
$$
h_{i}(x,t)=
\begin{cases}
~~ \frac{1}{3}\lambda u^2,&\text{$i$=0}, \\
~~ \frac{1}{3}\lambda u^2,&\text{$i$=1}, \\
-\frac{2}{3}\lambda u^2,&\text{$i$=2}.
\end{cases}
$$

\vspace{-0.5mm}
\section{Numerical results}\label{S3}
\vspace{-1mm}

Hereby, we give four computational experiments for Eq.~(\ref{Eq1}), the theoretical ones of which has already been presented in Ref.~\cite{21}. By comparing these two types of results, we can confirm the accuracy of the presented model. The velocity $u(x,t)$ is determined by incipient condition and the distribution function $f_i(x,t)$ is first initialized to be $f_{i}^{(0)}(x,t)$. In addition, the absolute error (AE) and global relative error (GRE) are presented to test the effectiveness of this model, which are defined as
\vspace{-3mm}
$$
\small {\rm AE}=|{u^ \dag}({x_i},t) - u({x_i},t)|, ~{\rm GRE} = \frac{{\sum\limits_i {|{u^ \dag}({x_i},t) - u({x_i},t)|} }}{{\sum\limits_i {|{u^ \dag}({x_i},t)|} }},
$$
where $u^ \dag (x_i,t)$ and $u(x_i,t)$ are the theoretical and computational solutions, respectively. In Ref.~\cite{21}, the nonlinear, dispersive and external-force coefficients of Eq.~(\ref{Eq1}) are given as
\vspace{-1mm}
\begin{equation}\label{Eq21}
\setlength\abovedisplayskip{0.5cm}
\setlength\belowdisplayskip{0.5cm}
\begin{aligned}
a(x,t)=\frac{x b_1(t)+b_2(t)}{p+ \int b_1(t)[q + \int m(t) dt] dt},
~b(x,t)=x b_1(t)+b_2(t),
~m(x,t)=m(t).
\end{aligned}
\end{equation}
The first-order theoretical solution of Eq.~(\ref{Eq1}) is given by~\cite{21}
\vspace{-1mm}
\begin{small}
\begin{equation}\label{Eq22}
\setlength\abovedisplayskip{0.8cm}
\setlength\belowdisplayskip{0.8cm}
\begin{aligned}
&u(x,t)=q+\int m(t) dt+L \tanh\Bigg[ w + \frac{L x}{2[p + \int b_1(t)[q + \int m(t) dt] dt]}
- \frac{L}{2} \int \frac{b_2(t)[q+\int m(t) dt]}{\big[p+\int b_1(t)[q+\int m(t) dt] dt \big]^2} dt\Bigg].&
\end{aligned}
\end{equation}
\end{small}

\vspace{-0.5mm}
\noindent\textbf{{\textbf Example 1.}}
In this simulation, we take the same parameters with Ref.~\cite{21}. Here $b_1(t)$, $b_2(t)$ and $m(t)$ are set to be 0, -2 and 0. Other parameters in our model are set as $w=6$, $p=-5$, $q=10$, $L=4$, $\Delta t=0.0001$ and $\Delta x=0.01 $, respectively. The simulation domain is immobilized to [0,40].
Then the theoretical solution is expressed as
\vspace{-1mm}
\begin{equation}\label{Eq23}
\begin{aligned}
u(x,t)=10+4\tanh(6+1.6t-0.4x).
\end{aligned}
\end{equation}
Here Eq.~(\ref{Eq1}) is restricted to the incipient condition
\vspace{-1mm}
\begin{equation}\label{Eq24}
\begin{aligned}
u(x,0)=10+4\tanh(6-0.4x),~0 \leq x\leq 40,
\end{aligned}
\end{equation}
\vspace{-2mm}
and the marginal conditions
\begin{equation}\label{Eq25}
\begin{aligned}
u(0,t)=14,~u(40,t)=6, ~t \geq 0.
\end{aligned}
\end{equation}
\vspace{-6mm}

\begin{figure}[!ht]
\renewcommand{\figurename}{Fig.}
\centering
\subfigure[]{\includegraphics[width=200 bp, height= 180bp]{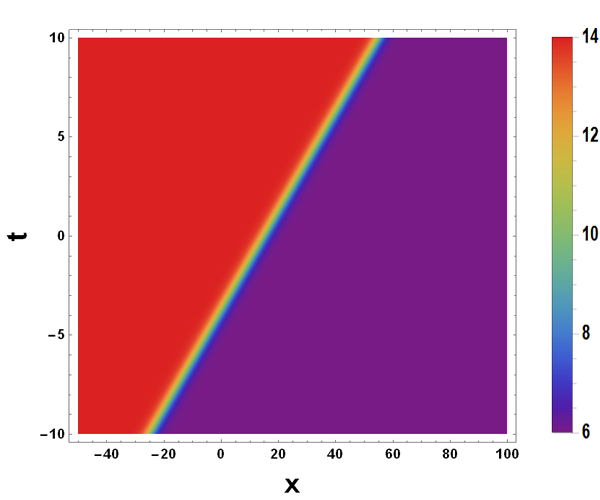}}
\qquad\qquad
\subfigure[]{\includegraphics[width=200 bp, height= 180bp]{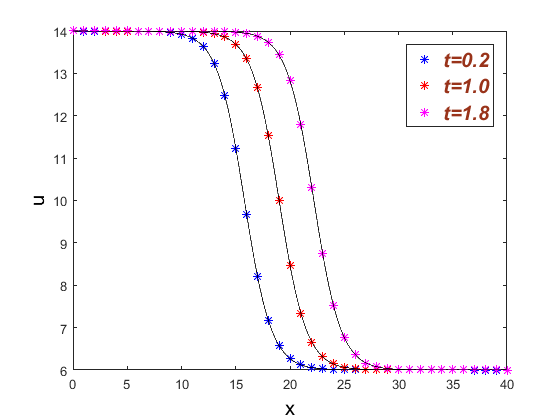}}
\vspace{-6mm}
\caption{\footnotesize (Color online) (a) Density plot of the theoretical solution~(\ref{Eq23}); (b) Computational results and the theoretical solution of Example 1 at $t=0.2, 1.0, 1.8$.
Discrete points and solid lines respectively denote the computational and theoretical solutions.}
\label{Fig1}
\end{figure}

\textbf{Fig.~\ref{Fig1}(a)} shows that the velocity of the soliton doesn't change with time~\cite{21}. \textbf{Fig.~\ref{Fig1}(b)} depicts that computational results align well with the theoretical ones. In detail, the GRE and relaxation time ($\tau$) of computational results are given in \textbf{Tab.~\ref{Tab1}}.
\vspace{-1mm}

\begin{table} [H] \scriptsize
\renewcommand{\tablename}{Tab.}
\caption{\footnotesize The GRE and relaxation time ($\tau$) for computational results of Example 1 }
\begin{tabular}{p{50mm}p{50mm}p{50mm}l}
\hline\noalign{\smallskip}\specialrule{0em}{1pt}{1pt}
$t$ & $t$=0.2 & $t$=1.0 & $t$=1.8  \\\specialrule{0em}{0.5pt}{0pt}
\noalign{\smallskip}\hline\noalign{\smallskip}\specialrule{0em}{0.5pt}{1pt}
GRE & 5.1826e-05 & 1.9123e-04 & 2.6970e-04\\\specialrule{0em}{0.5pt}{1pt}
$\tau$ & 2.5 & 2.5 & 2.5\\\specialrule{0em}{0.5pt}{0pt}
\noalign{\smallskip}\hline
\end{tabular}\label{Tab1}
\end{table}

\vspace{-1mm}
\noindent\textbf{{\textbf Example 2.}}
Compared with the situation in Example 1, we choose the same functions and parameters except that $b_2(t)=-2-t^2$.
Then the theoretical solution of this example is written as
\vspace{-2mm}
\begin{equation}\label{Eq26}
\begin{aligned}
u(x,t)=10+4\tanh(6+1.6t+\frac{4}{15}t^3-0.4x).
\end{aligned}
\end{equation}

\begin{figure}[!ht]
\renewcommand{\figurename}{Fig.}
\centering
\subfigure[]{\includegraphics[width=200 bp, height= 180bp]{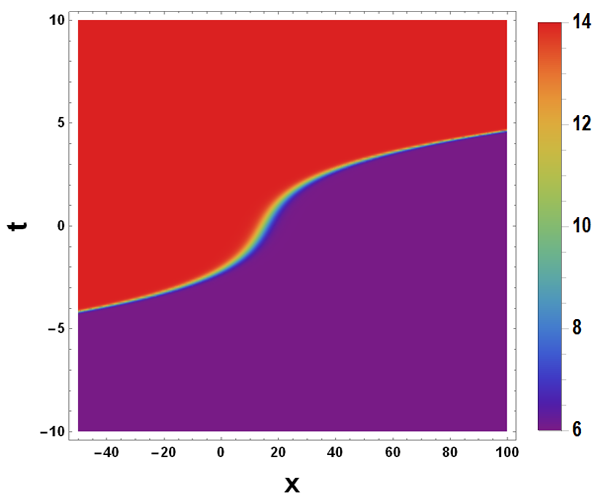}}
\qquad\qquad
\subfigure[]{\includegraphics[width=200 bp, height= 180bp]{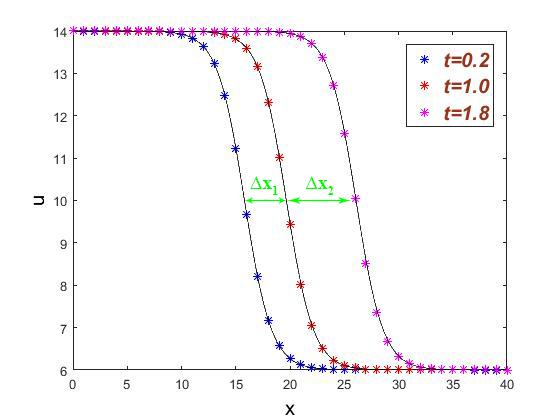}}
\vspace{-5mm}
\caption{\footnotesize (Color online) (a) Density plot of the theoretical solution~(\ref{Eq26}); (b) Computational results and the theoretical solution of Example 2 at $t=0.2, 1.0, 1.8$. The meanings of discrete points and solid lines are the same as {Fig.~\ref{Fig1}}.}
\label{Fig2}
\end{figure}

\vspace{-1.5mm}
The GRE and relaxation time ($\tau$) of computational results are given in \textbf{Tab.~\ref{Tab2}}.
In order to compare these two types of results explicitly, we give the detailed computational and theoretical results in \textbf{Tab.~\ref{Tab3}}. The AE shows that computational results are accurate enough.
\textbf{Fig.~\ref{Fig2}(a)} reveals that the dispersion changes the velocity of the soliton~\cite{21}. Here we use $\Delta x$ in \textbf{Fig.~\ref{Fig2}(b)} to represent the change of displacement during the same time $\Delta t=0.8$. We can find that computational results describe the effect of the dispersion on the motion of the soliton pretty well ($\Delta x_1<\Delta x_2$).

\vspace{-3mm}

\vspace{-8mm}

\vspace{-1mm}
\begin{table}[H]\scriptsize
\renewcommand{\tablename}{Tab.}
\caption{\footnotesize The GRE and  relaxation time ($\tau$) for computational results of Example 2 }
\begin{tabular}{p{50mm}p{50mm}p{50mm}l}
\hline\noalign{\smallskip}\specialrule{0em}{0.5pt}{1pt}
$t$ & $t$=0.2 & $t$=1.0 & $t$=1.8  \\ \specialrule{0em}{0.5pt}{0pt}
\noalign{\smallskip}\hline\noalign{\smallskip} \specialrule{0em}{0.5pt}{1pt}
GRE & 1.2185e-04 & 4.1204e-04 & 8.8659e-04\\ \specialrule{0em}{0.5pt}{1pt}
$\tau$ & 2.5 & 2.5225 & 2.86\\\specialrule{0em}{0.5pt}{0pt}
\noalign{\smallskip}\hline
\end{tabular}\label{Tab2}
\end{table}

\vspace{-5mm}
\begin{table}[H]\scriptsize
\renewcommand{\tablename}{Tab.}
\caption{\footnotesize The comparison between the theoretical and computational results of Example 2 at $t$=0.2 }
\begin{tabular}{p{50mm}p{50mm}p{50mm}l}
\hline\noalign{\smallskip}\specialrule{0em}{0.5pt}{1pt}
$x$ & theoretical results & computational results & AE  \\\specialrule{0em}{0.5pt}{0pt}
\noalign{\smallskip}\hline\noalign{\smallskip}\specialrule{0em}{0.5pt}{1pt}
4.0 & 13.999367 & 13.999370 & 3.0860e-06\\ \specialrule{0em}{0.5pt}{1pt}
8.0 & 13.984498 & 13.984510 & 1.2332e-05\\ \specialrule{0em}{0.5pt}{1pt}
12.0 & 13.636275 & 13.636260 & 1.4848e-05\\ \specialrule{0em}{0.5pt}{1pt}
16.0 & 9.6891613 & 9.678672 & 1.0489e-02\\ \specialrule{0em}{0.5pt}{1pt}
20.0 & 6.2696634 & 6.267866 & 1.7974e-03\\ \specialrule{0em}{0.5pt}{1pt}
24.0 & 6.0113594 & 6.011281 & 7.8366e-05\\ \specialrule{0em}{0.5pt}{1pt}
28.0 & 6.0004637 & 6.000460 & 3.6643e-06\\ \specialrule{0em}{0.5pt}{1pt}
32.0 & 6.0000189 & 6.000019 & 9.8970e-08\\ \specialrule{0em}{0.5pt}{1pt}
36.0 & 6.0000008 & 6.000001 & 2.2955e-07\\ \specialrule{0em}{0.5pt}{0pt}
\noalign{\smallskip}\hline
\end{tabular}\label{Tab3}
\end{table}

\vspace{-4mm}
\noindent\textbf{{\textbf Example 3.}}
In this example, we select the same functions and parameters with Example 1 except that $m(t)=\frac{1}{2}\sin(t+5)$.
Therefore, the expression of the theoretical solution is as follows
\vspace{-2mm}
\begin{equation}\label{Eq27}
\begin{aligned}
u(x,t)=10-\frac{1}{2}\cos(5+t)+4\tanh\big[6+1.6t-0.4x-\sin(5+t)\big].
\end{aligned}
\end{equation}
\vspace{-4mm}

\begin{figure}[!ht]
\renewcommand{\figurename}{Fig.}
\centering
\subfigure[]{\includegraphics[width=200 bp, height= 180bp]{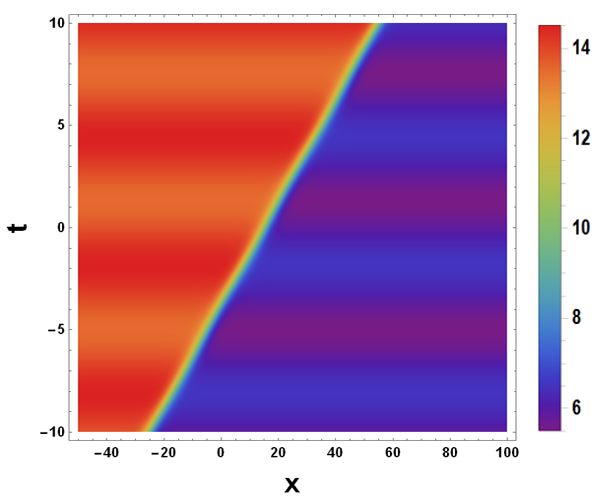}}
\qquad\qquad
\subfigure[]{\includegraphics[width=200 bp, height= 180bp]{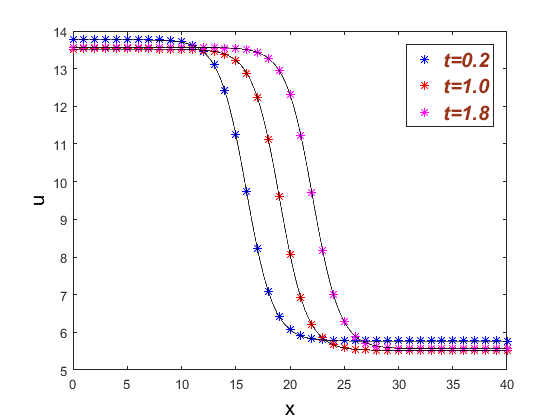}}
\vspace{-5mm}
\caption{\footnotesize (Color online) (a) Density plot of the theoretical solution~(\ref{Eq27}); (b) Computational results and the theoretical solution of Example 3 at $t=0.2, 1.0, 1.8$. The meanings of discrete points and solid lines are the same as {Fig.~\ref{Fig1}}.}
\label{Fig3}
\end{figure}

The GRE and relaxation time ($\tau$) of computational results are shown in \textbf{Tab.~\ref{Tab4}}.
Besides, we give the detailed results and AE in \textbf{Tab.~\ref{Tab5}}.
The external-force term changes the amplitude of the background of the soliton, which is shown in \textbf{Fig.~\ref{Fig3}(a)}~\cite{21}.
As depicted in \textbf{Fig.~\ref{Fig3}(b)}, computational results reflect such influence well.

\vspace{-5mm}
\begin{table}[H]\scriptsize
\renewcommand{\tablename}{Tab.}
\caption{\footnotesize The GRE and  relaxation time ($\tau$) for computational results of Example 3 }
\begin{tabular}{p{50mm}p{50mm}p{50mm}l}
\hline\noalign{\smallskip}\specialrule{0em}{0.5pt}{1pt}
$t$ & $t$=0.2 & $t$=1.0 & $t$=1.8  \\\specialrule{0em}{0.5pt}{0pt}
\noalign{\smallskip}\hline\noalign{\smallskip}\specialrule{0em}{0.5pt}{1pt}
GRE & 9.7518e-04 & 3.0205e-04 & 4.8431e-04\\\specialrule{0em}{0.5pt}{1pt}
$\tau$ & 2.5 & 2.5 & 2.5\\ \specialrule{0em}{0.5pt}{0pt}
\noalign{\smallskip}\hline
\end{tabular}\label{Tab4}
\end{table}
\vspace{-3mm}

\vspace{-2mm}
\begin{table}[H]\scriptsize
\renewcommand{\tablename}{Tab.}
\caption{\footnotesize The comparison between the theoretical and computational results of Example 3 at $t$=0.2}
\begin{tabular}{p{50mm}p{50mm}p{50mm}l}
\hline\noalign{\smallskip} \specialrule{0em}{0.5pt}{1pt}
$x$ & theoretical results & computational results & AE  \\ \specialrule{0em}{0.5pt}{1pt}
\noalign{\smallskip}\hline\noalign{\smallskip}  \specialrule{0em}{0.5pt}{1pt}
4.0 & 13.76519 & 13.774070 & 8.8803e-03\\ \specialrule{0em}{0.5pt}{1pt}
8.0 & 13.752222 & 13.761110 & 8.8881e-03\\ \specialrule{0em}{0.5pt}{1pt}
12.0 & 13.446754 & 13.455740 & 8.9865e-03\\ \specialrule{0em}{0.5pt}{1pt}
16.0 & 9.7284481 & 9.736932 & 8.4839e-03\\ \specialrule{0em}{0.5pt}{1pt}
20.0 & 6.0735016 & 6.082262 & 8.7604e-03\\ \specialrule{0em}{0.5pt}{1pt}
24.0 & 5.7787673 & 5.787643 & 8.8757e-03\\ \specialrule{0em}{0.5pt}{1pt}
28.0 & 5.7662734 & 5.775154 & 8.8806e-03\\ \specialrule{0em}{0.5pt}{1pt}
32.0 & 5.7657633 & 5.774644 & 8.8807e-03\\ \specialrule{0em}{0.5pt}{1pt}
36.0 & 5.7657425 & 5.774623 & 8.8805e-03\\ \specialrule{0em}{0.5pt}{0pt}
\noalign{\smallskip}\hline
\end{tabular}\label{Tab5}
\end{table}

\vspace{-8mm}

\vspace{-4mm}
\noindent\textbf{{\textbf Example 4.}}
Finally, we give the simulation for the theoretical solution with both the variable external-force and dispersive terms whose expression is presented by
\vspace{-2mm}
\begin{equation}\label{Eq28}
\begin{aligned}
&u(x,t)=10-\frac{1}{2}\cos(5+t)+4\tanh\big[6+1.6t+\frac{4}{15}t^3-0.4x+6\cos(5+t)+3t\sin(5+t)\big].&
\end{aligned}
\end{equation}

\begin{figure}[!ht]
\renewcommand{\figurename}{Fig.}
\centering
\subfigure[]{\includegraphics[width=200 bp, height= 180bp]{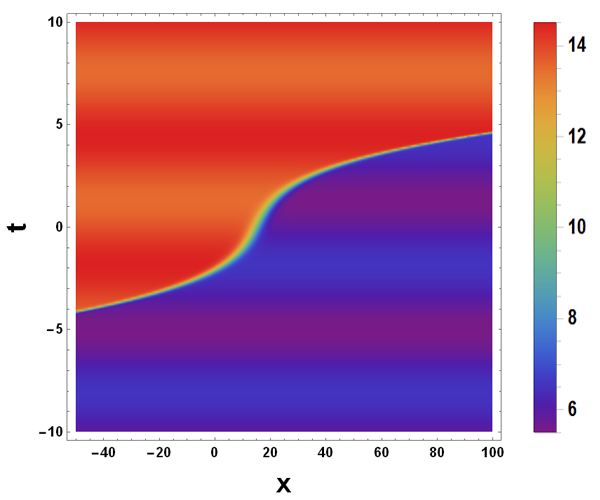}}
\qquad\qquad
\subfigure[]{\includegraphics[width=200 bp, height= 180bp]{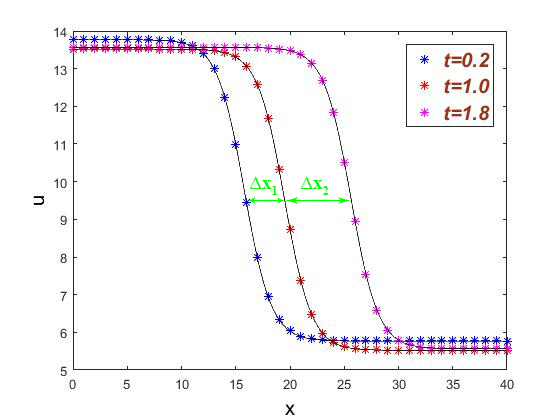}}
\vspace{-5mm}
\caption{\footnotesize (Color online) (a) Density plot of the theoretical solution~(\ref{Eq28}). (b) Computational results and theoretical solutions of Example 4 at $t=0.2, 1.0, 1.8$. The meanings of discrete points and solid lines are the same as {Fig.~\ref{Fig1}}.}
\label{Fig4}
\end{figure}

\vspace{-1mm}
\textbf{Fig.~\ref{Fig4}(a)} shows that these two variable coefficients change both the amplitude and velocity ($\Delta x_1<\Delta x_2$)~\cite{21}, and the definition of $\Delta x$ is the same as Example 2. As shown in \textbf{Fig.~\ref{Fig4}(b)}, this situation combines two characteristics in Example 2 and 3. The GRE and relaxation time ($\tau$) of computational results are provided in \textbf{Tab.~\ref{Tab6}}.

\vspace{-1mm}

\vspace{-5mm}
\begin{table}[H]\scriptsize
\renewcommand{\tablename}{Tab.}
\caption{\footnotesize The GRE and relaxation time ($\tau$) for computational results of Example 4 }
\begin{tabular}{p{50mm}p{50mm}p{50mm}l}
\hline\noalign{\smallskip} \specialrule{0em}{0.5pt}{1pt}
$t$ & $t$=0.2 & $t$=1.0 & $t$=1.8  \\ \specialrule{0em}{0.5pt}{0pt}
\noalign{\smallskip}\hline\noalign{\smallskip} \specialrule{0em}{0.5pt}{1pt}
GRE & 9.7555e-04 & 2.6925e-04 & 5.3918e-04\\ \specialrule{0em}{0.5pt}{1pt}
$\tau$ & 2.5324 & 3.4604 & 5.6684\\ \specialrule{0em}{0.5pt}{0pt}
\noalign{\smallskip}\hline
\end{tabular}\label{Tab6}
\end{table}

\vspace{-6mm}
\section{Conclusions}\label{S4}
\vspace{-3mm}

In summary, the lattice Boltzmann model has been constructed for Eq.~(\ref{Eq1}) by choosing the equilibrium distribution and compensatory functions. Eq.~(\ref{Eq1}) has been recovered correctly through the Chapman-Enskog expansion. Four numerical examples have been studied when we take some appropriate parameters and they agree well with the theoretical results. Through our work, the dynamic characteristics affected by the dispersive and external-force coefficients, such as the velocity and amplitude, have been simulated effectively. Based on the error analysis, it has been found that the LBM could be applied to the variable-coefficient NLEEs with the external-force term.

\bibliography{abanov-bibliography}

\end{spacing}
\end{document}